\documentclass[12pt]{amsart}
\usepackage{comment,amsmath,amstext,amsfonts,amssymb,graphicx,framed,dsfont,enumitem,xcolor,tikz}

\usepackage{mathrsfs}
\usepackage{float}

\theoremstyle{definition}
\newtheorem{prop}{Proposition}

\newtheorem{lemma}[prop]{Lemma}

\newtheorem{theo}[prop]{Theorem}

\DeclareSymbolFontAlphabet{\mathrsfs}{rsfs}

\newcommand{\h}{\mathcal{H}}
\newcommand{\degh}{\deg_{\h}}
\newcommand{\deh}{\mathrm{d}_{\h}}
\newcommand{\rh}{\mathrm{r}(\h)}
\newcommand{\arh}{\mathrm{ar}(\h)}
\newcommand{\qh}{\mathrm{q}(\h)}
\newcommand{\Dunh}{\Delta(\h)}
\newcommand{\Ddeuxh}{\Delta([\h]_2)}

\begin{document}

\title{About Berge-F\"uredi's conjecture on the \\ chromatic index of hypergraphs}
\author{Alain Bretto}
\address{NormandieUnicaen, GREYC CNRS-UMR 6072, Caen, France}
\email{alain.bretto@info.unicaen.fr}

\author{Alain Faisant}
\author{François Hennecart}
\address{
Universit\'e Jean Monnet, ICJ UMR5208, CNRS, Ecole Centrale de Lyon, INSA Lyon, Universite Claude Bernard Lyon 1,
42023 Saint-Etienne, France}
\email{faisant@univ-st-etienne.fr, francois.hennecart@univ-st-etienne.fr}

\subjclass[2010]{05C15}
\keywords{Hypergraphs, hyperedge coloring, chromatic index}

\begin{abstract} 
We show that the chromatic index of a hypergraph $\h$ satisfies
{Berge-F\"uredi} conjectured bound $\qh\le  \Ddeuxh+1$ under certain hypotheses on the antirank $\arh$ or on the maximum degree $\Dunh$. This provides sharp information in connection with 
{Erd\H{o}s-Faber-Lov\'asz} Conjecture which deals with the coloring of a family of cliques that intersect pairwise in at most one vertex.
\end{abstract}

\maketitle

\section{\bf Introduction}

Berge-F\"uredi's conjecture (cf. \cite{BER1989,BER1997,FUR1986} and {Meyniel}, unpublished communication)  extends to linear hypergraphs {Vizing}'s theorem on the chromatic index of graphs (see \cite[Theorem 7.2.2]{BFH2022}). It asserts that the chromatic index 
of any linear hypergraph $\h=(V, E)$  without loop satisfies 
\begin{equation}\label{eq1}
\qh\le\Ddeuxh+1.
\end{equation}
This bound was conjectured first around 1986. Since then, to our knowledge, only few advances have been made. Recently in \cite{ZHANGSKINNER2020}, it is shown that \eqref{eq1} holds whenever $\Ddeuxh\le 5$.
In the present paper we obtain new conditional results where the antirank $\arh$ and the maximum degree $\Dunh$ are involved together with the maximum degree $\Ddeuxh$ of the $2$-section of $\h$.  

If $\h$ is linear then $\Ddeuxh\le |V|-1$. Hence bound \eqref{eq1} implies {Erd\H os-Faber-Lov\'asz} Conjecture asserting that $\qh\le |V|$ for any linear hypergraph $\h$.

\section{\bf Notations and statements}

Let $\h=(V,E)$ be a hypergraph. 
For $x\in V$ and $e\in E$,
we let 
$$
\degh(x)=|\{a\in E\,:\, x\in a\}|,\quad
\deh(e)=|\{a\in E\smallsetminus\{e\}\,:\, a\cap e\ne\varnothing\}|,
$$
be the \emph{degree of vertex $x$} and the \emph{degree of hyperedge $e$} in $\h$, respectively. 
We then define respectively the \emph{rank}, the \emph{antirank}, the \emph{maximum degree} and the \emph{minimum degree} of $\h$ by
$$
\rh=\max_{e\in E}|e|,\quad 
\arh=\min_{e\in E}|e|,\quad 
\Dunh=\max_{x\in V}\degh(x),\quad
\delta(\h)=\min_{x\in V}\degh(x).
$$ 
A \emph{loop} in a hypergraph $\mathcal{H}$ is a hyperedge $e$ such that $\vert e\vert=1$. 
\\
The \emph{$2$-section} of $\h=(V,E)$ is the loopless multigraph $[\h]_2=(V,A)$ with $k$ many  edges between $x,y\in V$, $x\ne y$, if there are $k$ hyperedges in $\h$ containing both $x$ and~$y$.  We have
\begin{equation}\label{eqn0}
\Ddeuxh\ge (\arh-1)\Dunh.
\end{equation}
The \emph{chromatic index} of $\h$, denoted by
 $\qh$, is  the least integer $q$ such 
that there exists a map $c:E\longrightarrow\{1,2,\dots,q\}$ satisfying
$$
\forall e,e'\in E,\ \left(c(e)=c(e')\implies e\cap e' =\varnothing \text{ or } e=e'\right).
$$
A   hypergraph   $\h=(V,E)$ is said \emph{linear} if for any pair of distinct  hyperedges $e, e'\in E$, we have $\vert e\cap e'\vert\leq 1$. In this case $[\h]_2$ reduces to a simple graph. 

\medskip
Recently it has been established that if $\h$ is a linear hypergraph 
satisfying 
\begin{equation}\label{arh}
\arh\ge \sqrt{|V|},
\end{equation} 
then $\qh\le |V|$
(cf. \cite{WANG2024}). This upper bound has been improved in \cite{BFH2024} where it is shown that under the same hypotheses on $\h$, then \eqref{eq1} actually holds. In this note we show that condition \eqref{arh} can be weakened 
by proving the following sharper result.

\begin{theo}\label{thm1}
Let $\h=(V,E)$ be a loopless hypergraph  such that $\arh \ge \sqrt{\Ddeuxh+1}$.
Then $\qh\le \Ddeuxh+1$.
\end{theo}

Note that we do not assume in the above statement that $\h$ is linear.
We will propose an approach supported by the notion of critical hypergraph
(see also \S\ref{rk81} for additional comments).

We will prove the
following key result which has its own interest.
$\h$ is said \emph{$k$-uniform} (resp. \emph{$d$-regular}) if 
$\forall e\in E$, $|e|=k$ (resp. $\forall x\in V,\ \degh(x)=d$).

\begin{theo}\label{thm2}
Let $k\ge2$ and $\h=(V,E)$ be a  $k$-uniform and $(k+1)$-regular linear hypergraph.
Then $\qh\le \Ddeuxh+1$.
\end{theo}

The following statement will be  derived from Theorem \ref{thm1}. 

\begin{theo}\label{thm3}
Let $\mathcal{H}=(V, E)$ be a loopless hypergraph such that
\begin{equation}\label{eqcor1}
\Dunh\leq \sqrt{\Ddeuxh+1}+1.
\end{equation}
Then
$\qh\leq \Ddeuxh+1.$
\end{theo}

\section{\bf Proof of Theorem \ref{thm1}}

In \cite{BFH2024} we exploited the notion of critical hyperedge.
A hyperedge $e$ is said \emph{critical} if  
$\mathrm{q}(\h\smallsetminus e)=\qh-1$. 
Hypergraph $\h$ is a \emph{critical hypergraph} if all its hyperedges are critical.
Since the addition of $e$ to $\h\smallsetminus e$ needs a new color, $e$ must intersect every color of $\h\smallsetminus e$. So we have the following.

\begin{lemma}[cf. Lemma 3.2 of \cite{BFH2024}]
\label{lemkey}
Let $\h=(V,E)$ be a hypergraph without loop and $e\in E$ be a critical hyperedge. Then
$\qh-1\le \mathrm{d}_{\h}(e)$.  
\end{lemma}

We may assume that $\h$ is connected and critical
since $\mathrm{ar}(\h')\ge \arh$ and $\Delta([{\h}']_2)\le \Delta([{\h}]_2)$ for any  partial hypergraph ${\h}'$ of ${\h}$. 
Let $e\in E$ be of minimum cardinality $k:=\arh$. 
We have 
\begin{equation}\label{eqdhe}
\mathrm{d}_{\h}(e)=\sum_{\substack{a\in E\smallsetminus\{e\}\\a\cap e\ne\varnothing}}1\le \sum_{a\in E\smallsetminus\{e\}}|a\cap e|
=\sum_{x\in V}\sum_{\substack{a\in E\smallsetminus\{e\}\\a\ni x}}1
=\sum_{x\in e}(\mathrm{deg}_{{\h}}(x)-1),
\end{equation}
and equality occurs when $\h$ is linear.
Then by Lemma \ref{lemkey} we get
$$
\mathrm{q}({\h})-1=\mathrm{q}({\h}\smallsetminus e)\le \deh(e)
\le k(\Delta({\h})-1)=(k-1)\Delta({\h})-(k-\Delta({\h})),
$$
so 
\begin{equation}\label{eq25}
\qh\le \Ddeuxh+1 -(k-\Dunh).
\end{equation}
Our hypothesis implies $\Ddeuxh\le k^2-1$. We distinguish two cases.

\begin{itemize}
\item If $\Ddeuxh<k^2-1$ we get by \eqref{eqn0} that
$\Dunh<k+1$ thus $\Dunh\le k$, yielding the required upper bound
for $\qh$ 
by \eqref{eq25}.

\item Assume that $\Ddeuxh=k^2-1$. We thus have to prove that $\qh\le k^2$.\\
If there exists $e\in E$ such that $|e|=k$ and $\deh(e)\le k^2-1$ then
$\qh\le \deh(e)+1\le k^2$ by Lemma \ref{lemkey} and we are done. 
Otherwise 
\begin{equation}\label{eqn1}
\deh(e)\ge k^2\quad\text{ for any hyperedge $e$ of cardinality $k$.}
\end{equation}
We shall see that $\h$ is a very special hypergraph.\\[0.2em]
\indent -- Let $e\in E$ be of cardinality $k$.  By \eqref{eqn0}  we have
$\Dunh\le k+1$, thus by \eqref{eqn1} and \eqref{eqdhe} we get
$$
k^2\le \deh(e)\le\sum_{x\in e}(\mathrm{deg}_{{\h}}(x)-1)\le k(\Dunh-1)
\le k(k+1-1)=k^2.
$$
We infer
\begin{equation}\label{eqn2}
\deh(e)=k^2\quad \text{and}\quad
\mathrm{deg}_{{\h}}(x)-1=k\quad  \text{for any $x\in e$.}
\end{equation}
\indent -- Let $a\in E$ such that $a\cap e\ne \varnothing$. Let $x\in a\cap e$.
We have
\begin{align*}
k^2-1=\Ddeuxh\ge \mathrm{deg}_{[\h]_2}(x)&=\sum_{\substack{b\in E\\b\ni x}}(|b|-1)\\
&=|a|-1+\sum_{\substack{b\in E\smallsetminus\{a\}\\b\ni x}}(|b|-1)\\
&\ge |a|-1+(\degh(x)-1)(k-1)\\
&=
 |a|-1+k(k-1)\\
&\ge k^2-1,
\end{align*}
giving $|a|=k$. \\[0.2em]
\indent -- Since $\h$ is connected we infer that 
$|a|=k$ for any $a\in E$. Hence $\h$ is $k$-uniform and $(k+1)$-regular.
By \eqref{eqdhe} and \eqref{eqn2}, we have $|a\cap e|=1$ whenever $a\cap e\ne\varnothing$. Thus $\h$ is linear. The proof will be completed when Theorem \ref{thm2} will be established.
\end{itemize}

\section{\bf Proof of Theorem \ref{thm2}}

By assumption we have 
\begin{equation}\label{eq21}
k|E|=(k+1)|V|,\quad 
\Ddeuxh=k^2-1\quad\text{and} \quad  \forall e\in E,\ d_{\h}(e)=k^2.
\end{equation}
Since $\h$ is linear 
its \emph{line-graph} $\mathrm{L}({\h})$ 
is the simple graph whose vertices are hyperedges of ${\h}$ and
such that there is an edge between $e,e'\in E$, $e\ne e'$, 
if $e\cap e'\ne \varnothing$. We have 
$\mathrm{deg}_{\mathrm{L}({\h})}(e)=\mathrm{d}_{\h}(e)$ for any $e\in E$.

By \eqref{eq21},  each vertex
$e$ of $\mathrm{L}(\h)$ has degree 
$\mathrm{deg}_{\mathrm{L}(\h)}(e)=k^2$. Hence 
$\mathrm{L}(\h)$ is $k^2$-regular. Since $\qh$ coincides with 
$\chi(\mathrm{L}(\h))$, namely the chromatic number of the simple graph 
$\mathrm{L}(\h)$, we infer from Brooks' theorem 
(see \cite[Theorem 7.1.9, p. 238]{BFH2022}) that either 
$$
\qh\le \Delta(\mathrm{L}(\h))=k^2=\Ddeuxh+1$$ 
or $\mathrm{L}(\h)$
is complete or an odd length cycle. 
If $\mathrm{L}(\h)$ were a complete graph, we would have $|E|=k^2+1$
implying by \eqref{eq21} that  $(k+1)$ divides $k(k^2+1)$, a contradiction since $k\ge2$. If $\mathrm{L}(\h)$ were a cycle, we would have
$d_{\h}(e)=\mathrm{deg}_{\mathrm{L}(\h)}(e)=2$ for any $e\in E$, a contradiction.  This ends the proof of Theorem \ref{thm2}.

\section{\bf Proof of Theorem \ref{thm3}}

Let ${\h}'=(V,E')$ be a critical partial hypergraph of ${\h}$ such that 
$\mathrm{q}({\h}')=\mathrm{q}({\h})$. 

\begin{itemize}
\item If $\mathrm{ar}({\h}')\ge\sqrt{\Delta([\h']_2)+1}$ we infer from Theorem \ref{thm1} that 
$$
\mathrm{q}({\h})=\mathrm{q}({\h}')\le \Delta([\h']_2)+1\le \Ddeuxh+1.
$$

\item 
If $\mathrm{ar}({\h}')<\sqrt{\Delta([\h']_2)+1}$ then we take $e\in E'$ such that
$|e|=\mathrm{ar}({\h}')$.  By Lemma \ref{lemkey} and using \eqref{eqdhe} with $\h'$ in place of $\h$, we get
\begin{align*}
\mathrm{q}({\h})-1=\mathrm{q}({\h}')-1&\le \mathrm{d}_{\h'}(e)\\
&\le \sum_{x\in e}(\mathrm{deg}_{{\h}'}(x)-1)\\
&\le |e|(\Delta({\h}')-1)\\
&=
\mathrm{ar}({\h}')(\Delta({\h'})-1)\\
&< \sqrt{\Delta([\h']_2)+1}(\Delta({\h'})-1)\\
&\le \sqrt{\Ddeuxh+1}(\Delta({\h})-1)\\
&\le \Ddeuxh+1
\end{align*}
by \eqref{eqcor1}.
 Hence $\qh\leq  \Ddeuxh+1$ as wanted.
\end{itemize}

\section{\bf Remarks and questions}

\subsection{}\label{rk81} We should notice that when 
$\arh>\sqrt{\Ddeuxh+1}$ then the desired 
bound $\qh\le 
\Ddeuxh+1$ 
can be  equally
deduced from greedy algorithm for coloring hyperedges (see for instance \cite[page 586]{DVOR2000}) which yields
$$
\qh\le \max_{\rh\ge k\ge \arh} k\left\lfloor\frac{\Ddeuxh}{k-1}-1\right\rfloor+1.
$$
We use here a different approach which allows us to solve in addition the equality case
$\arh=\sqrt{\Ddeuxh+1}$.

\subsection{}\label{rk82}
When the maximum degree $\Dunh$ satisfies
 $\Dunh > \sqrt{\Ddeuxh+1}+1$, Theorem~\ref{thm3} is no longer usable. However at the cost of an additional condition on the rank $\rh=\max_{e\in E}|e|$, we can notice that \eqref{eq1} remains valid for linear hypergraphs $\h$.
Let $\delta_2:=\sqrt{\Ddeuxh}$ and 
assume that 
 \begin{align*}
\rh&=\delta_2-u,\  u>0,\\
\Dunh &= \delta_2 +1+v,\ v>0.
\end{align*}
Then 
$$
\qh= \chi(\mathrm{L}(\h))\leq \Delta(\mathrm{L}(\h))+1
\leq \rh(\Dunh-1)+1=(\delta_2-u)(\delta_2+v)+1.
$$
Thus \eqref{eq1} holds whenever $(v-u)\delta_2\le uv$, that is
$\frac1u\le \frac1v+\frac1{\delta_2}$.

\subsection{}
If $\h$ is supposed linear, in the \emph{extremal case} $\Ddeuxh=|V|-1$ , $\arh=|V|^{1/2}$ in Theorem \ref{thm1} the authors obtained in \cite{BFH2024}
a much better upper bound for $\qh$, namely 
$\qh\le \frac12(\Ddeuxh+1+\sqrt{\Ddeuxh+1})+1$.
 How the upper bound for $\qh$ can be improved  under the sole condition~ $\arh = \sqrt{\Ddeuxh+1}$?

\subsection{}
Theorem \ref{thm2} deals with $k$-uniform and $d$-regular linear hypergraphs
when $d=k+1$ and provides the required bound \eqref{eq1}. 
For general $k,d$ satisfying $k\le d-2$,  the line-graph $\mathrm{L}(\h)$ is $k(d-1)$ regular and Brooks' theorem  provides the non sufficiently sharp $\qh=\chi(\mathrm{L}(\h))\le \Delta(\mathrm{L}(\h))=k(d-1)$  while $\Ddeuxh=d(k-1)$.

\subsection{}
More generally let $\h$ be a $k$-uniform linear hypergraph with $k\ge2$ and 
set $\delta_2:=\sqrt{\Ddeuxh}$. Then \eqref{eq1} holds if at least one of the following condition is satisfied.
\begin{enumerate}
\item $k=2$: Vizing's theorem  (see \cite[Theorem 7.2.2, p. 249]{BFH2022}); 
\item $k\ge\sqrt{\Ddeuxh+1}$: Theorem \ref{thm1};
\item $\Ddeuxh=k^2$ which implies $\Dunh\le \sqrt{\Ddeuxh}+1$: Theorem \ref{thm3};
\item $3\le k <\delta_2$ and $\frac1{\delta_2-k}\le \frac1{\Dunh-\delta_2-1}+\frac{1}{\delta_2}$: see \S\ref{rk82}.
\end{enumerate}
In order to solve Berge-F\"uredi's conjecture for all linear uniform hypergraphs,
it thus remains to consider 
the case where 
$$
3\le k<\delta_2<\Dunh-1  \quad\text{and}\quad
\frac1{\delta_2-k} > \frac1{\Dunh-\delta_2-1}+\frac{1}{\delta_2}.
$$
We may likewise ask the following question:
is Berge-F\"uredi's conjecture in all its generality implied by its validity for linear uniform hypergraphs ?


\end{document}